# Numerical Study of Cosserat Fluid-Structure Interaction in a Monolithic Eulerian Framework

Nazim Hussain Hajano[a], Muhammad Sabeel Khan[b,*], Lisheng Liu[c]

[a]Department of Engineering Structure and Mechanics, Wuhan University of Technology, Wuhan, People's Republic of China, 430070.
[b]Department of Mathematics, Sukkur Institute of Business Administraiton University, 65200 Sukkur, Pakistan
[c]Department of Engineering Structure and Mechanics, Wuhan University of Technology, Wuhan, People's Republic of China, 430070.

*Corresponding author: m.sabeel@iba-suk.edu.pk

**Abstract**

We propose a monolithic Eulerian variational formulation in non-classical sense of continuum description for the analysis of micro-viscosity parameters at micro-structural level. In this respect, Cosesrat fluid-strucutre interaction CFSI phenomena is taken into account by considering micro-rotational degrees of freedom of the particles. The governing equations and variational formulation for CFSI problem are presented. The space and time domains are discretized by the finite element method and semi-implicit scheme. The model is implemented and evaluated using FreeFem++. The present model is analyzed by computing a well known benchmark problem FLUSTRUK-FSI-3*, where the flow around a flag attached with cylinder is considered for numerical test. The obtained results suggest that the amplitude of oscillations and micro-rotaional viscosity $\mu_r$ varies inversely and micro-viscosity parameter $\lambda_1$ exhibit a significant effect on micro-rotation field as compare to velocity field.

**Graphical Abstract**

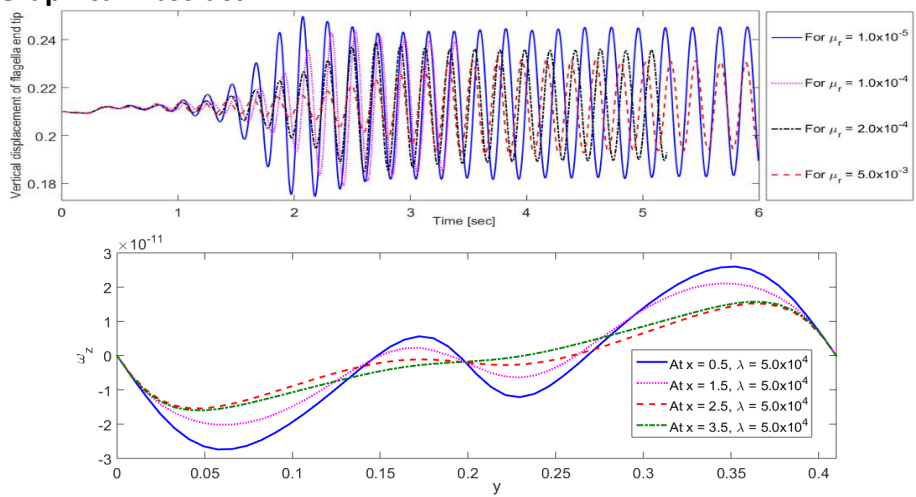

# 1 INTRODUCTION

A monolithic approach finds origin in (Hron and Turek 2006) solve fluid-structure interaction (FSI) problems with a single variational equation for the whole FSI system (Dunne 2006; Heil et al. 2008), or more recently (Pironneau 2016; Hecht and Pironneau 2017; Chiang et al. 2017; Pironneau 2018; Murea 2019). Among these, a monolithic Eulerian approach (Pironneau 2016; Hecht and Pironneau 2017) is similar to the fully Eulerian formulation (Dunne and Rannacher 2006; Dunne 2006) presents FSI system equations in terms of displacement, velocity and pressure as usual but in an Eulerian framework. The major difference between two formulations is that in (Dunne and Rannacher 2006; Dunne 2006) prefer to work with velocities in the fluid domain and displacement in the solid domain, while monolithic Eulerian approach (Pironneau 2016; Hecht and Pironneau 2017) proposes to work with velocities everywhere in the problem domain. The fully Eulerian formulation replaces the well-established ALE, formulates fluid and solid sub domains in Eulerian framework, respectively. More work devoted to fully Eulerian formulation can be seen in (Richter and Wick 2010; Rannacher and Richter2011; Richter 2013; Wick 2013).

All numerical approaches exsist in literature deals with classical continuum description of FSI problems. Here in this study, a monolithic approach is presented in the framework of non-classical Cosserat continuum description for solving and analyszing micro-structural FSI phenomena.

Classical continuum mechanics consider continuum as a simple point-continua with points having three displacement (dofs) and a symmetric Cauchy stress tensor that characterizes the response of a material to the displacement. Such classical models may not be sufficient for the description of non-classical physical phenomena, where microstructural effects are observed most in high strain gradients regions. The Cosserat continuum theory (Eugene and François Cosserat 1909) is one of the most prominent theory to model non-classical physical phenomena. Further, this concept was applied to describe fluids with microstructures by (Condiff and Dahler 1964; Eringen 1964 &1966) and the mathematical details, with some of its applications, are presented in (Lukaszewicz 1999). In non-classical continuum theory, the response of the material to the displacement and micro-rotation is characterized by a nonsymmetric Cauchy stress tensor and couple stresss tensor, respectively.

In the present study, a monolithic Eulerian approach is employed to analyze the effect of micro-viscosity parameters on the Cosserat fluid-structure interaction CFSI problem which is not yet studied in the literature. This benchmark problem was first studied by (Schafer and Turek 1996), and later by (Turek and Hron 2006; Dunne and Rannacher 2006; Hecht and Pironneau 2017), respectively. All the above considerations were in the case of Classical continuum framework but here in this study a non-classical continuum framework is utilized to study the flow around the cylinder in particular Cosserat continuum. The algorithmic description is presented and implemented using FreeFem++ (Hecht 2012).

This paper is organized as follows. In section 2, we describe the Continuum description, notations, the constitutive relations and the derivation of the governing equations from the conservation laws for the Cosserat model. In Section 3, we present the variational formulation in monolithic Eulerian framework for the CFSI problem. In Section 4-5, the time and spatial discretization is presented using semi-implicit scheme and the finite element method, respectively. Numerical tests and results are discussed in detailed in section 6 and section 7. Future developments and conclusion of the study is addressed in section 8.

# 2 CONTINUUM MECHANICS

## 2.1 Notations

Consider $\Omega^t$ time dependent computational domain comprises of the fluid region $\Omega_f^t$ and solid region $\Omega_s^t$ such that $\bar{\Omega}^t = \bar{\Omega}_f^t \cup \bar{\Omega}_s^t, \Omega_f^t \cap \Omega_s^t = \emptyset, \forall t$. The notations $\Sigma^t = \bar{\Omega}_f^t \cup \bar{\Omega}_s^t$ and $\partial \Omega^t$ denotes fluid-structure interface and the boundary of computational domain $\Omega^t$, respectively. The part of $\partial\Omega^t$ on which either the structure is clamped or on which there is 'no slip condition' on the fluid, that part is denoted by $\Gamma$ assumed to be independent of time $t$. At initial time the fluid domain $\Omega_f^0$ and the solid domain $\Omega_s^0$ are prescribed.

The following are considered standard notations as in (Ciarlet 1988; Marsden and Hughes 1993; Bath 1996; Tallec and Mouro 2001; Hron and Turek 2006). We use bold characters for vectors and tensors, with some exceptions, like $x, x^0 \in \mathbb{R}^d, d = 2$ or $3$.

- $\mathbf{X} : \Omega^0 \times (0, T) \mapsto \Omega^t : \mathbf{X}(x^0, t)$, the Lagrangian position at time $t$ of $x^0$.
- $\mathbf{d} = \mathbf{X}(x^0, t) - x^0$, the displacement.
- $\mathbf{u} = \partial_t \mathbf{X}$, the velocity of the deformation.
- $\mathbf{F}_{ji} = \partial_{x_i^0} \mathbf{X}_j$, the transposed gradient of deformation.
- $J = \det \mathbf{F}$, the Jacobian of the deformation.

Let $\rho_f$ and $\sigma_f$ be the density and stress tensor in the fluid, as for the solid structure with $\rho_s, \sigma_s$, respectively. We define

- $\rho(x, t) = \rho_f \mathbf{1}_{\Omega_f^t}(x, t) + \rho_s \mathbf{1}_{\Omega_s^t}(x, t)$, the density.
- $\boldsymbol{\sigma}(x, t) = \boldsymbol{\sigma}_f \mathbf{1}_{\Omega_f^t}(x, t) + \boldsymbol{\sigma}_s \mathbf{1}_{\Omega_s^t}(x, t)$, the stress tensor.

Where $\mathbf{1}_{\Omega_f^t}$ is the set function indicator defined as $\mathbf{1}_\Omega(x) = \begin{cases} 1 & if\ x \in \Omega \\ 0 & otherwise \end{cases}$.

Unless specified otherwise, all spatial derivatives are with respect to $x \in \Omega^t$ and not with respect to $x^0 \in \Omega^0$. If $\phi$ is a function of $x = \mathbf{X}(x^0, t)$ where $x^0 \in \Omega^0$, then

$$\nabla_{x^0} \phi = \left[\partial_{x_i^0} \phi\right] = \left[\partial_{x_i^0} \mathbf{X}_j \partial_{x_j} \phi\right] = \mathbf{F}^T \nabla \phi. \tag{1}$$

When $\mathbf{X}$ is one-to-one and invertible, the relation between $\mathbf{F}$ and $\mathbf{d}$ can be seen as function of $(x, t)$ instead of $(x^0, t)$. Mathematically

$$\mathbf{F}^T = \nabla_{x^0} \mathbf{X} = \nabla_{x^0}(\mathbf{d} + x^0) = \nabla_{x^0} \mathbf{d} + \mathbf{I} = \mathbf{F}^T \nabla \mathbf{d} + \mathbf{I} \quad \Rightarrow \quad \mathbf{F} = (\mathbf{I} - \nabla \mathbf{d})^{-T}. \tag{2}$$

The time convective derivative of $\phi$ becomes

$$D_t \phi := \frac{\partial}{\partial t} \phi(\mathbf{X}(x^0, t), t) = \partial_t \phi(x, t) + \mathbf{u} \cdot \nabla \phi(x, t). \tag{3}$$

Finally, we introduce the deformation tensor and the micro-rotaion strain tensor, respectively as

$$\mathbf{Du} = \nabla \mathbf{u} + \nabla \mathbf{u}^T, \tag{4}$$

$$\kappa = \nabla \boldsymbol{\omega}. \tag{5}$$

**2.2 Conservation Laws**

Conservation of mass, Linear momentum and angular momentum for fluid and solid medium takes form

$$\frac{d}{dt}(J\rho) = 0, \tag{6}$$

$$\rho D_t \mathbf{u} = \nabla \cdot \boldsymbol{\sigma}_{s,f} + \mathbf{f}, \tag{7}$$

$$\rho I \dot{\boldsymbol{\omega}} = \nabla \cdot \mathbf{C}_f + \boldsymbol{\varepsilon} : \boldsymbol{\sigma}_f + \mathbf{g}. \tag{8}$$

Where $\mathbf{u}$ and $\boldsymbol{\omega}$ are velocity and microrotation fields, respectively. The body force density, the body couple density and the microinertia coefficient are denoted by $\mathbf{f}$, $\mathbf{g}$ and $I$ respectively. Finally, incompressibility implies $J=1$ and $\rho = \rho^0$.

## 2.3 Constitutive Equations

Two-dimensional hyperelatic incompressible Mooney-Rivlin material and incompressible viscous Cosserat fluid is considered.

- For an incompressible hyper-elastic material the constitutive description is stated as

$$\boldsymbol{\sigma}_s = -p_s \mathbf{I} + \rho_s \frac{\partial \Psi}{\partial \mathbf{F}} \mathbf{F}^T, \tag{9}$$

such that the Helmoltz potential $\Psi$ in of the Mooney-Rivilin material is given by

$$\Psi(\mathbf{F}) = c_1 tr(\mathbf{F}^T \mathbf{F}) + c_2 \left\{ tr(\mathbf{F}^T \mathbf{F})^2 - tr^2(\mathbf{F}^T \mathbf{F}) \right\}. \tag{10}$$

Where the constants $c_1$ and $c_2$ are empically determined (Ciarlet 1988).

Consequently, $\qquad \partial_{\mathbf{F}} \Psi \mathbf{F}^T = (2c_1 - 4c_2 c) \mathbf{B} + 4c_2 \mathbf{B}^2. \tag{11}$

Taking $\mathbf{B} = \mathbf{FF}^T = \left[ (\mathbf{I} - \nabla \mathbf{d})(\mathbf{I} - \nabla \mathbf{d})^T \right]^{-1}, b = \det \mathbf{B}, c = tr\mathbf{B}, \mathbf{B}^2 = c\mathbf{B} - b\mathbf{I},$ and applying Calay-Hamilton theorem, equation (10) becomes

$$\partial_{\mathbf{F}} \Psi \mathbf{F}^T = 2c_1 (\mathbf{Dd} - \nabla \mathbf{d} \nabla \mathbf{d}^T) + \alpha' \mathbf{I}, \tag{12}$$

Where $\alpha' = 2c_1 tr(\mathbf{FF}^T) - (2c_1 + 4c_2) \det(\mathbf{FF}^T)$ is some scalar function of material parameters $c_1$ and $c_2$.

- For incompressible viscous Cosserat fluid the constitutive relations are described as

$$\boldsymbol{\sigma}_f = -p_f \mathbf{I} + \mu (\nabla \mathbf{u} + \nabla \mathbf{u}^T) + \mu_r (\nabla \mathbf{u} - \nabla \mathbf{u}^T) - 2\mu_r \boldsymbol{\varepsilon} \cdot \boldsymbol{\omega}, \tag{13}$$

$$\mathbf{C}_f = \alpha (tr\boldsymbol{\kappa}) \mathbf{I} + \beta (\boldsymbol{\kappa} + \boldsymbol{\kappa}^T) + \gamma (\boldsymbol{\kappa} - \boldsymbol{\kappa}^T). \tag{14}$$

Where $\boldsymbol{\sigma}_f$, $\mathbf{C}_f$ and $\mathbf{I}$ are the nonsymmetric stress tensor, the couple stress tensor and the identity tensor, respectively. The pressure field and coefficient of dynamic viscosity are denoted by $p$ and $\mu$, respectively. The coefficients of microviscosity are represented as $\mu_r$, $\alpha$, $\beta$ and $\gamma$. The Levi Civita tensor is represented by $\varepsilon$.

## 2.4 Derivation of Governing Equations for Cosserat Fluid

Subjected to certain prescribed boundary conditions according to the description of the physical problem and taking into consideration the constitutive equations, the governing conservation equations (6) - (8) for Cosserat fluid as in the book (Lukaszewicz 1999) leads to

$$\nabla \cdot \mathbf{u} = 0, \tag{15}$$

$$\rho\left(\frac{\partial \mathbf{u}}{\partial t} + \mathbf{u}\cdot\nabla\mathbf{u}\right) = -\nabla p + (\mu+\mu_r)\Delta\mathbf{u} + 2\mu_r \nabla\times\boldsymbol{\omega} + \mathbf{f}, \tag{16}$$

$$\rho\left(I\frac{\partial \boldsymbol{\omega}}{\partial t} + I\mathbf{u}\cdot\nabla\boldsymbol{\omega}\right) = \lambda_1\Delta\boldsymbol{\omega} + \lambda_2\nabla\nabla\cdot\boldsymbol{\omega} - 4\mu_r\boldsymbol{\omega} + 2\mu_r(\nabla\times\mathbf{u}) + \mathbf{g}. \tag{17}$$

Where $\lambda_1 = \beta+\gamma$ and $\lambda_2 = \alpha+\beta-\gamma$ are material parametres related to micro-viscosity.

## 3 Monolithic Variational Formulation

This variational formulation considers homogeneous boundary conditions on $\Gamma \subset \partial\Omega^t$ and homogeneous Neumann conditions on $\partial\Omega^t \setminus \Gamma$. In monolithic Eulerian variational formuation the final CFSI in 2-dimensions is:

Given $\Omega_f^0$, $\Omega_s^0$ and $\mathbf{d}$, $\mathbf{u}$ at $t=0$, find $\mathbf{u}, \boldsymbol{\omega}, p, \mathbf{d}, \Omega_f^t, \Omega_s^t$ with $\mathbf{u}|_\Gamma = 0$ and $\boldsymbol{\omega}|_\Gamma = 0$ such that

$$\int_{\Omega_f\cup\Omega_s} \rho D_t\mathbf{u}\cdot\tilde{\mathbf{u}} - p\nabla\cdot\tilde{\mathbf{u}} - \tilde{p}\nabla\cdot\mathbf{u} + (\mu+\mu_r)\mathbf{Du}:\mathbf{D\tilde{u}} - 2\mu_r(\nabla\times\boldsymbol{\omega})\cdot\tilde{\mathbf{u}}\, d\Omega$$
$$+ \int_{\Omega_s} c_3(\mathbf{Dd} - \nabla\mathbf{d}\nabla\mathbf{d}^T):\mathbf{D\tilde{u}}\, d\Omega_s = \int_{\Omega_f\cup\Omega_s} \mathbf{f}\cdot\tilde{\mathbf{u}}\, d\Omega, \tag{18}$$

$$\int_{\Omega_f\cup\Omega_s} \rho D_t I\boldsymbol{\omega}\cdot\tilde{\boldsymbol{\omega}} + \lambda_1\nabla\boldsymbol{\omega}:\nabla\tilde{\boldsymbol{\omega}} - \lambda_2\nabla\cdot\boldsymbol{\omega}\cdot\tilde{\boldsymbol{\omega}} + 4\mu_r\boldsymbol{\omega}\cdot\tilde{\boldsymbol{\omega}} - 2\mu_r\nabla\times\mathbf{u}\cdot\tilde{\boldsymbol{\omega}}\, d\Omega = \int_{\Omega_f\cup\Omega_s} \mathbf{g}\cdot\tilde{\boldsymbol{\omega}}\, d\Omega, \tag{19}$$

for all $(\tilde{\mathbf{u}}, \tilde{\boldsymbol{\omega}}, \tilde{p})$ with $\tilde{\mathbf{u}}|_\Gamma = 0$ and $\tilde{\boldsymbol{\omega}}|_\Gamma = 0$, where $\Omega_f^t$ and $\Omega_s^t$ are defined incrementally by $D_t\mathbf{d} = \mathbf{u}$

$$\text{and} \qquad \frac{dX}{d\tau} = \mathbf{u}(X(\tau),\tau), \quad X(t)\in\Omega_r^t \Rightarrow X(\tau)\in\Omega_r^\tau \quad \forall\tau\in(0,T), \quad r=s,f. \tag{20}$$

The relation in (20) defines $\Omega_f^t$ and $\Omega_s^t$ forward in time. Above the notations $\mathbf{B}:\mathbf{C} = tr(\mathbf{B}^T\mathbf{C})$ and $c_3 = \rho^s c_1$ are used.

## 4 TIME DISCRETIZATION

A semi-implicit scheme is used to discretize CFSI equations based on its monolithic Eulerian variational formulation (18-19). In this respect, let $t\in[0,T]$ be the simulation time where $T$ is total time. Discretize the interval $[0,T]$ into equall sub intervals each of the length $\delta t = \frac{T}{N}$ such that $t = n\delta t$ where $n = 0, 1, \cdots, N$.

Let $\mathbf{d}^{n+1} = \mathbf{d}^n + \delta t\mathbf{u}^{n+1}$. Hence

$$\mathbf{Dd} - \nabla\mathbf{d}\nabla\mathbf{d}^T \approx \mathbf{Dd}^n - \nabla\mathbf{d}^n\mathbf{d}^{n^T} + \delta t(\mathbf{Du}^{n+1} - \nabla\mathbf{u}^{n+1}\mathbf{d}^{n^T} - \nabla\mathbf{d}^n\nabla\mathbf{u}^{n+1^T}) + o(\delta t). \tag{21}$$

Now, if $X^n$ is a first order approximation of $X(t^{n+1} - \delta t)$ defined by $\dot{X} = \mathbf{u}(X(\tau),\tau)$, $X(t^{n+1}) = x$, where $X(t^{n+1}) = x$ such that $X^n(x) = x - \delta t \mathbf{u}^n(x)$, then a first order in time approximation for CFSI problem (18-19) reads:

Find $\mathbf{u}^{n+1} \in \mathbf{H}_0^1(\Omega^{n+1})$, $\boldsymbol{\omega}^{n+1} \in \mathbf{H}_0^1(\Omega^{n+1})$, $p^{n+1} \in L_0^2(\Omega^{n+1})$, $\Omega_f^{n+1}$ and $\Omega_s^{n+1}$ such that with $\mathbf{u}^{n+1}|_\Gamma = 0$, $\boldsymbol{\omega}^{n+1}|_\Gamma = 0$ and $\Omega^{n+1} = \Omega_f^{n+1} \cup \Omega_s^{n+1}$, $\forall\, \tilde{\mathbf{u}} \in \mathbf{H}_0^1(\Omega^{n+1})$, $\tilde{\boldsymbol{\omega}} \in \mathbf{H}_0^1(\Omega^{n+1})$, $\tilde{p} \in L_0^2(\Omega^{n+1})$ with $\tilde{\mathbf{u}}|_\Gamma = 0$ and $\tilde{\boldsymbol{\omega}}|_\Gamma = 0$, the following holds

$$\int_{\Omega_f^n \cup \Omega_s^n}\left[\left(\rho^n \frac{\mathbf{u}^{n+1}-\mathbf{u}^n \circ X^n}{\delta t}\right)\cdot\tilde{\mathbf{u}} - p^{n+1}\nabla\cdot\tilde{\mathbf{u}} - \tilde{p}\nabla\cdot\mathbf{u}^{n+1} + (\mu+\mu_r)\mathbf{D}\mathbf{u}^{n+1}:\mathbf{D}\tilde{\mathbf{u}} - 2\mu_r(\nabla\times\boldsymbol{\omega}^{n+1})\cdot\tilde{\mathbf{u}}\right]d\Omega^n$$
$$+\int_{\Omega_s^n} c_3\left[\mathbf{D}\mathbf{d}^n - \nabla\mathbf{d}^{n^T}\nabla\mathbf{d}^n + \delta t\left(\mathbf{D}\mathbf{u}^{n+1} - \nabla\mathbf{u}^{n+1}\nabla\mathbf{d}^{n^T} - \nabla\mathbf{d}^n\nabla\mathbf{u}^{n+1^T}\right)\right]:\mathbf{D}\tilde{\mathbf{u}}\, d\Omega_s^n = \int_{\Omega_f^n \cup \Omega_s^n}\mathbf{f}\cdot\tilde{\mathbf{u}}\,d\Omega^n,\tag{22}$$

and

$$\int_{\Omega_f^n \cup \Omega_s^n}\left\{\rho^n I^n\left(\frac{\boldsymbol{\omega}^{n+1}-\boldsymbol{\omega}^n\circ \tilde{X}^n}{\delta t}\right)\cdot\tilde{\boldsymbol{\omega}} + \lambda_1 \nabla\boldsymbol{\omega}^{n+1}:\nabla\tilde{\boldsymbol{\omega}} - \lambda_2\nabla(\nabla\cdot\boldsymbol{\omega}^{n+1})\cdot\tilde{\boldsymbol{\omega}} + 4\mu_r\boldsymbol{\omega}^{n+1}\cdot\tilde{\boldsymbol{\omega}} - 2\mu_r(\nabla\times\mathbf{u}^{n+1})\cdot\tilde{\boldsymbol{\omega}}\right\}d\Omega^n$$
$$= \int_{\Omega_f^n \cup \Omega_s^n}\mathbf{g}\cdot\tilde{\boldsymbol{\omega}}\,d\Omega^n.\tag{23}$$

Now, update $\mathbf{d}$ by $\mathbf{d}^{n+1} = \mathbf{d}^n\circ X^n + \delta t\mathbf{u}^{n+1}$, and $\Omega_r^n$ by $\Omega_r^{n+1} = \{x + \delta t\mathbf{u}^{n+1}(x) : x \in \Omega_r^n\}$, where $r = s, f$.

## 5 SPATIAL DISCRETIZATION

The finite element method is used to discretize space domain. Let $V_h$ and $W_h$ represents the finite element functional spaces for the velocities, displacements and micro-rotational velocities, respectively and $Q_h$ be the functional space for pressure field. Let $\mathfrak{I}_h^0$ be a triangulation of the initial domain $\Omega^0$ with quadratic elements for displacements, translational, micro-rotational velocities and linear elements for pressure field. Given that that the pressure is different in fluid domain and structural domain because of the discontinuity of pressure at the fluid-structure interface $\Sigma$; therefore, the functional space $Q_h$ is space of piecewise linear functions on the triangulation and is continuous in $\Omega_r^{n+1}, r = s, f$. A small penalization parameter $\zeta \ll 1$ needs to be added to impose uniqueness of the pressure when one desire to use direct linear solver. The discrete variational formulation of CFSI, thus reads

Find $\mathbf{u}_h^{n+1}, \boldsymbol{\omega}_h^{n+1}, p_h^{n+1} : \forall\, \tilde{\mathbf{u}}_h \in V_{0h}, \tilde{\boldsymbol{\omega}}_h \in \tilde{W}_{0h}, \tilde{p}_h \in Q_h$ with $V_{0h}|_\Gamma = 0$ and $W_{0h}|_\Gamma = 0$ are subspaces of $V_h$ and $W_h$, such that

$$\int_{\Omega_f^n \cup \Omega_s^n}\left[\left(\rho^n \frac{\mathbf{u}_h^{n+1}-\mathbf{u}_h^n \circ X^n}{\delta t}\right)\cdot\tilde{\mathbf{u}}_h - p_h^{n+1}\nabla\cdot\tilde{\mathbf{u}}_h - \tilde{p}_h\nabla\cdot\mathbf{u}_h^{n+1} + (\mu+\mu_r)\mathbf{D}\mathbf{u}_h^{n+1}:\mathbf{D}\tilde{\mathbf{u}}_h - 2\mu_r(\nabla\times\boldsymbol{\omega}_h^{n+1})\cdot\tilde{\mathbf{u}}_h\right]\Omega^n$$
$$+\int_{\Omega_s^n} c_3\left[\mathbf{D}\mathbf{d}_h^n - \nabla\mathbf{d}_h^{n^T}\nabla\mathbf{d}_h^n + \delta t\left(\mathbf{D}\mathbf{u}_h^{n+1} - \nabla\mathbf{u}_h^{n+1}\nabla\mathbf{d}_h^{n^T} - \nabla\mathbf{d}_h^n\nabla\mathbf{u}_h^{n+1^T}\right)\right]:\mathbf{D}\tilde{\mathbf{u}}_h\, d\Omega_s^n + \int_{\Omega_f^n \cup \Omega_s^n}\zeta p_h\tilde{p}_h\,d\Omega^n = \int_{\Omega_f^n \cup \Omega_s^n}\mathbf{f}\cdot\tilde{\mathbf{u}}_h\,d\Omega^n,\tag{24}$$

and

$$\int_{\Omega_f^n \cup \Omega_s^n} \left\{ \rho^n I^n \left( \frac{\boldsymbol{\omega}_h^{n+1} - \boldsymbol{\omega}_h^n \circ \tilde{X}^n}{\delta t} \right) \cdot \tilde{\boldsymbol{\omega}}_h + \lambda_1 \nabla \boldsymbol{\omega}_h^{n+1} : \nabla \tilde{\boldsymbol{\omega}}_h - \lambda_2 \nabla \left( \nabla \cdot \boldsymbol{\omega}_h^{n+1} \right) \cdot \tilde{\boldsymbol{\omega}}_h + 4\mu_r \boldsymbol{\omega}_h^{n+1} \cdot \tilde{\boldsymbol{\omega}}_h - 2\mu_r \left( \nabla \times \mathbf{u}_h^{n+1} \right) \cdot \tilde{\boldsymbol{\omega}}_h \right\} d\Omega^n$$

$$= \int_{\Omega_f^n \cup \Omega_s^n} \mathbf{g} \cdot \tilde{\boldsymbol{\omega}}_h d\Omega^n. \quad (25)$$

To update the triangulation at each vertex (say $q_i^n$) of the triangle $T_h \in \Im_h^n$, the vertex is moved to a new position by

$$q_i^{n+1} := q_i^n + \delta t \mathbf{u}_h^{n+1}.$$

By denoting $\mathbf{d}_i^n := \mathbf{d}^n(q_i)$, it can be seen that

$$\mathbf{d}^n \circ X^n (q_i^{n+1}) = \mathbf{d}^n (q_i^n + \delta t \mathbf{u}_h^{n+1} - \delta t \mathbf{u}_h^{n+1}) = \mathbf{d}^n (q_i^n).$$

This implies that the displacement vector of vertices $\mathbf{d}_h^n$ can be copied to $\mathbf{d}_i^{n+1}$ plus with addition of $\delta t \mathbf{u}_h^{n+1}(q_i^n)$ in order to obtain $\mathbf{d}_h^{n+1}$, i.e.,

$$\mathbf{d}_h^{n+1} = \mathbf{d}_h^n \circ X^n + \delta t \mathbf{u}_h^{n+1}(q_i^n) = \mathbf{d}_h^n + \delta t \mathbf{u}_h^{n+1}(q_i^n).$$

Moreover, the fluid domain mesh is moved by $\tilde{\mathbf{u}}$ which is a solution of the Laplace problem $-\Delta \tilde{\mathbf{u}} = 0$ $\forall \tilde{\mathbf{u}} \in V_{0h}$, subjected to $\tilde{\mathbf{u}}|_\Sigma = \mathbf{u}$ where $\Sigma$ is Cosserat fluid structure interface and $\tilde{\mathbf{u}} = 0$ at the boundaries $\Gamma^f \cup \Gamma^s \setminus \Sigma$. Moving the vertices of each triangle $T_h \in \Im_h^n$ by the above procedure gives a new triangulation $\Im_h^{n+1}$.

# 6 NUMERICAL TESTS

The monolithic Eulerian formulation is used to analyze the effect of micro-viscosity parameters at micro-structural level for CFSI problem. The spatial discretization is performed by using Lagrangian triangular finite elements with quadratic elements for displacements, translational, micro-rotational velocities and linear elements for pressure field. The publically available tool Freefem++ (Hecht 2012) has been used to implement the algorithms.

## 6.1 Incompressible hyperelastic flag attached to a cylinder

The present model is analyzed by computing a well known benchmark problem FLUSTUK-FSI-3* (Hecht and Pironneau 2017). This benchmark problem was first studied by (Schafer and Turek 1996), and later by (Turek and Hron 2006; Dunne and Rannacher 2006; Hecht and Pironneau 2017), respectively. The description of the model problem for the present study is described as (Hecht and Pironneau 2017). An incompressible hyperelastic Mooney-Rivlin material, like rectangular flag of size $[0,l] \times [0,h]$, is attached at the back of a hard fixed cylinder in the computational rectangular domain $[0,L]\times[0,H]$. The fluid flow enters from left inlet and is free to leave on the right outlet.

## 6.2 Mathematical Description

We detail the following numerical values for geometry, boundary and intial conditions for the present problem as suggested in (Hecht and Pironneau 2017).

*Geometry:* The point $0.2, 0.2$ is the center of a cylinder with radius $r = 0.05$; the dimension of the rectangular flag and the fluid computational domain are considered in computation, respectively as $l = 0.35$, $h = 0.02$, $L = 2.5$ and $H = 0.41$ which set cylinder slightly below the symmetry line of the computational domain. The control point $A(t)$ is fixed at the trailing edge of structure with $A(0) = (0.6, 0.2)$.

*Boundary and Initial Conditions:* Top and bottom boundaries satisfies the 'no-slip' condition. A parabolic velocity profile is prescribed at the left inlet,

$$\mathbf{u}_f\left(0,y\right) = \bar{\mathbf{U}}\left(\frac{6y\,H-y}{H^2}\right),$$

where $\bar{\mathbf{U}}$ is mean inflow velocity with flux $\bar{\mathbf{U}}H$ and $\bar{\mathbf{U}} = 2$. The zero-stress $\sigma \cdot n = 0$ is employed at the right outlet using do-nothing approach. Furthermore, the density and the reduced kinematic viscosity of the fluid takes values $\rho_f = 10^3 kgm^{-3}$ and $\nu_f = \frac{\mu}{\rho_f} = 10^{-3} m^2 s^{-1}$. For solid structure we consider $\rho_s = \rho_f$, $c_1 = 10^6 kgm^{-1}s^{-1}$ and no external force. Initially, all flow velocities and structure displacements are zero.

## 7 RESULTS AND DISCUSSION

In this section, we discuss the effects of micro-viscosity parameters $\mu_r$ and $\lambda_1$ at micro-structural level for the present CFSI problem using obtained results. During numerical simulation flow starts oscillating around $t \sim 2$ and develops a Karman vortex street. The Figure 1, shows the effect of micro-rotational viscosity $\mu_r$ on vertical displacement of flagella end tip (control poin $A$) as a function of time for CFSI problem. The amplitude of the oscillations and micro-rotaional viscosity varies inversely, which specifies that fluctuations of the control poin $A$ almost disappear by dominating micro-rotational viscosity $\mu_r$ on classical viscosity $\mu$ of the fluid.

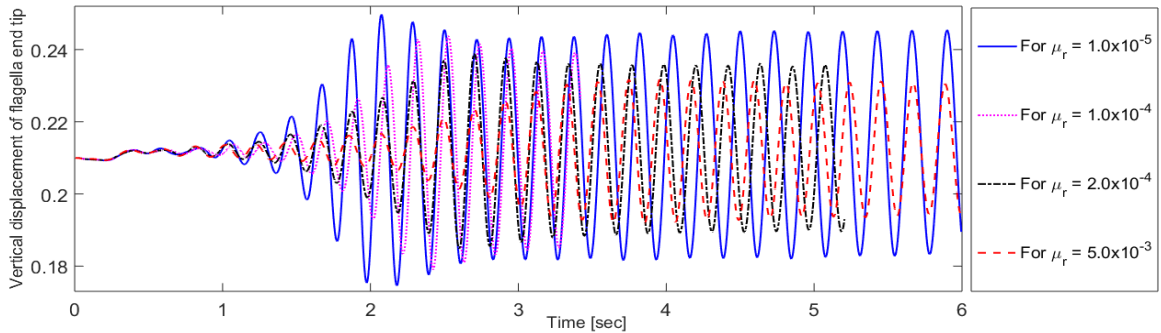

*Figure 1: Effect of micro-rotional viscosity on vertical displacement of flagella end tip different time.*

The fluid parameter $\lambda_1 = \lambda$ in our notation combines the shear spin and rotational spin viscosities, respectively, shows significant effect on micro-rotation field $\omega_z$. For a constant value of parameter $\lambda_1$, fluid particles experiences large rotational effect near the control point $A$ and negligible while moving beyond it in flow direction; consequently vanishes at the boundaries of the computational domain as shown in Figure 2. It is also observed that micro-rotational effect decreases with increasing value of fluid parameter $\lambda_1$ as displayed in Figure 3.

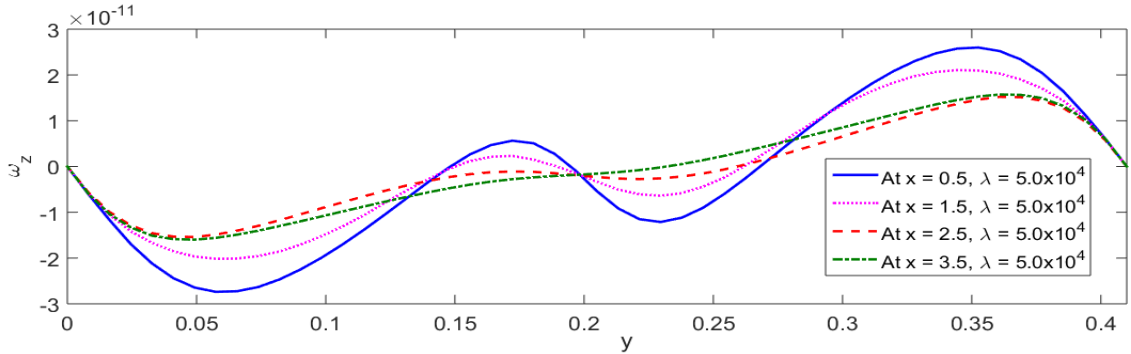

**Figure 2**: *Micro-rotational velocity at different horizontal positions with fixed $\lambda_1$.*

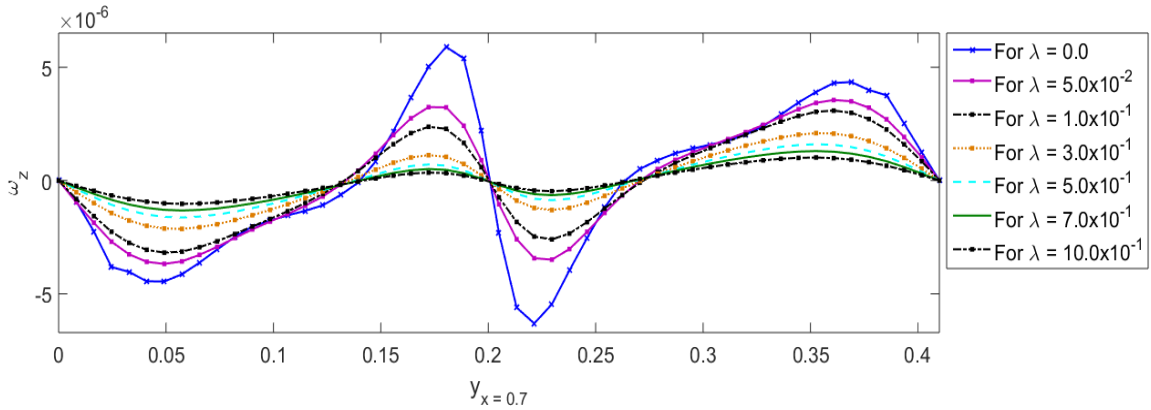

**Figure 3**: *The effects of $\lambda_1$ on micro-rotational velocity near the control point $A$ at $x = 0.7$.*

Consequently, on the bases of results shown in Figure 2 and Figure 3, it is further observed that

- Fluid particles near the bottom boundary wall of computational domain preserve clockwise rotation and anticlockwise rotation near the line $y = 0.2$; and follows the same pattern above the line $y = 0.2$.

- The rotational effect becomes minimum near the line $y = 0.2$ and almost vanishes on it in the computational domain.

Figure 4, suggests that the effect of $\lambda_1$ on the horizontal velocity component of fluid particles is negligible at all points in computational domain with symmetrical behavior. Moreover, fluid particles attains parabolic graph by increasing the $\lambda_1$ beyond the control point $A$, which results in maximum velocity magnitude on the symmetry line.

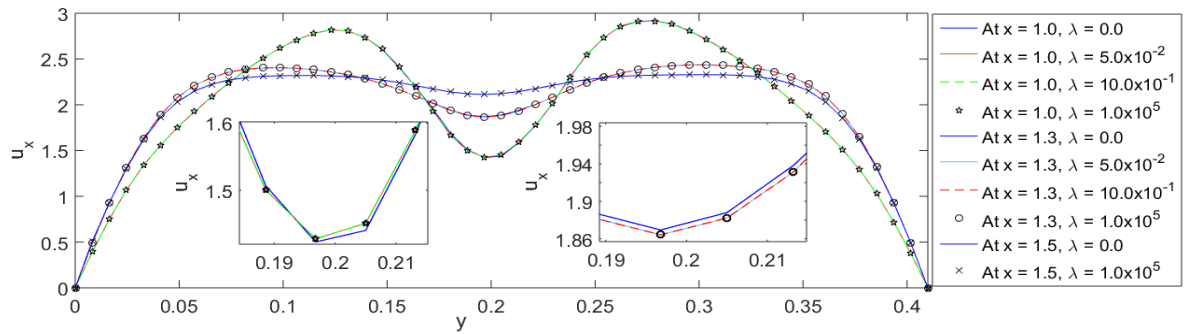

**Figure 4**: The effects $\lambda_1$ on $x$- component of velocity at different horizontal positions.

Furthermore, it is observed that micro-viscosity parameter $\lambda_1$ effects vertical component of velocity in same way but does not follow the particular pattern at different fixed horizontal positions beyond the control point $A$ as shown in Figure 5. Finally, the results disussed above demonstrate that micro-viscosity parameter $\lambda_1$ exhibit a significant effect on micro-rotation fieled as compare to velocity field. All results are obtaind by using the Reynolds number $Re = 200$ and mean velocity $\bar{U} = 2$,

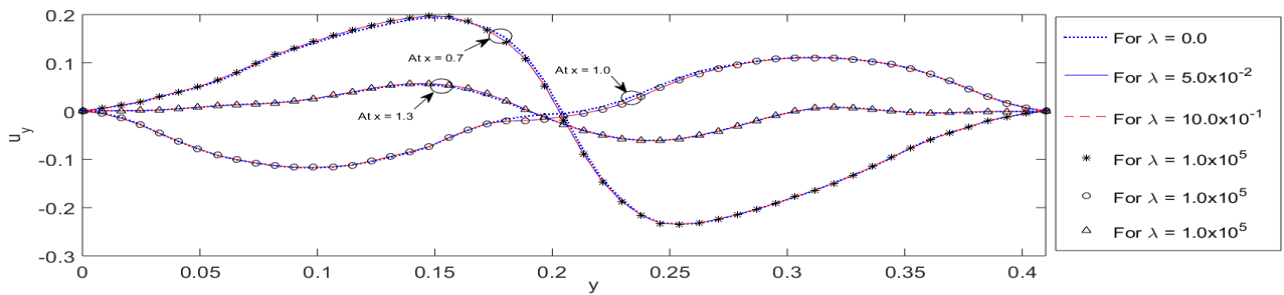

**Figure 5**: The effects of $\lambda_1$ on $y$ - component of velocity at different horizontal positions.

The Figure 6-7, represents some color snapshots of micro-rotational velocity profile and its contour plots during simulation at different time.

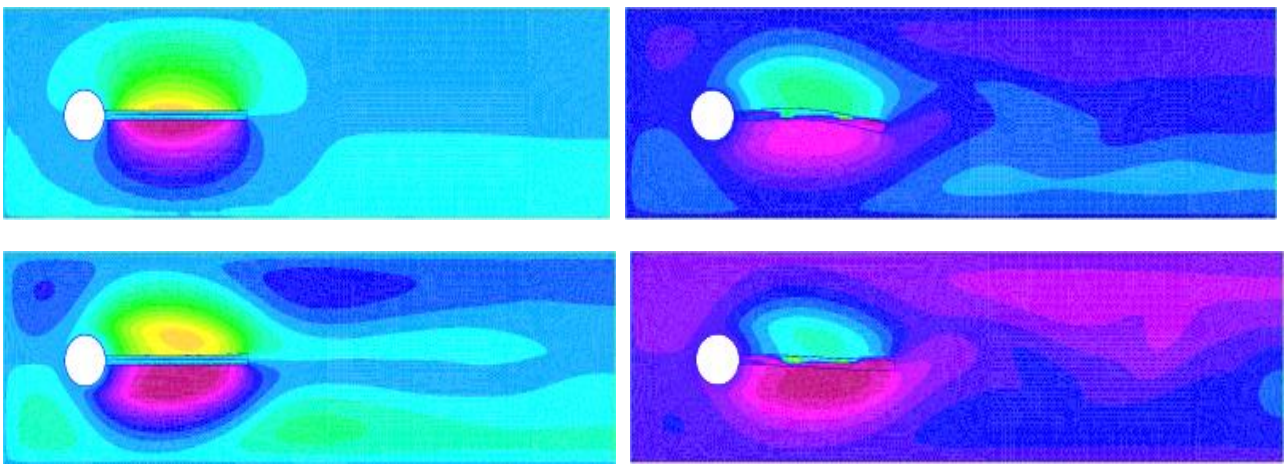

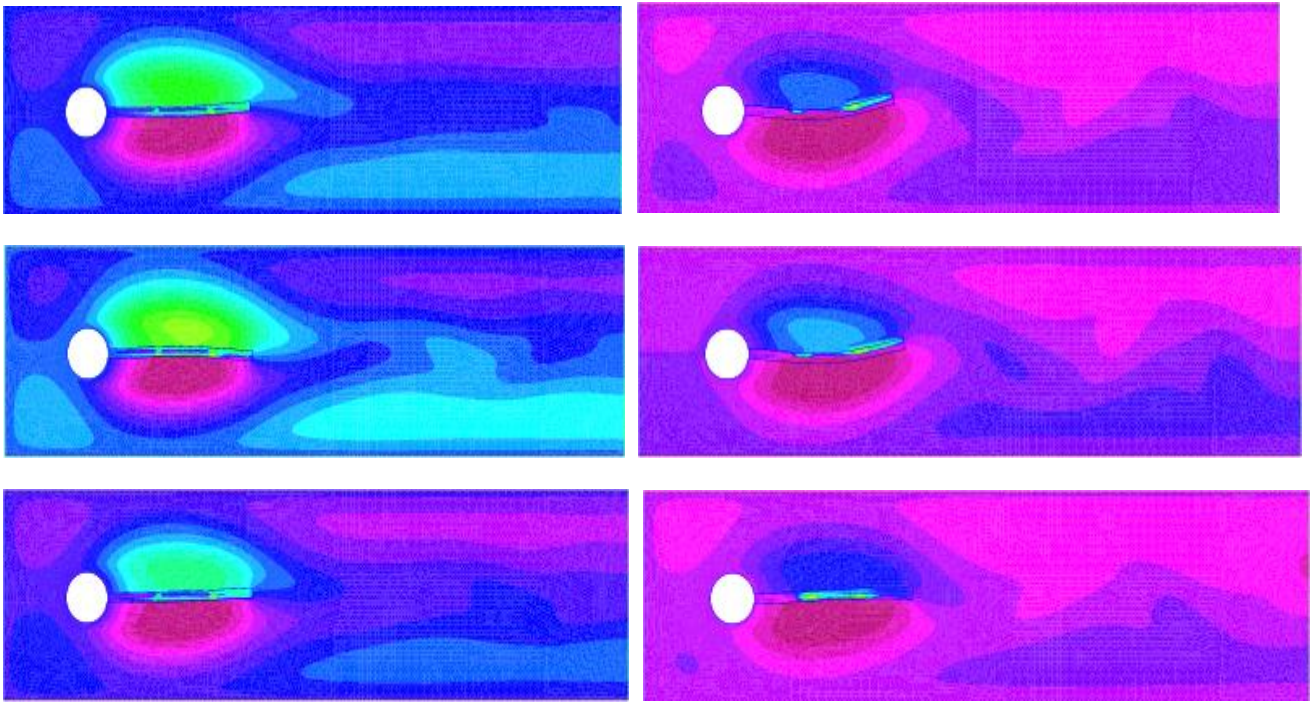

**Figure 6**: Micro-rotational velocity profile $w_h$ plots (*column* 1: $t = 0.005, 0.5, 1.8, 1.935, 1.975$ & *column* 2: $t = 2.3, 2.44, 2.65, 2.9, 3.0$).

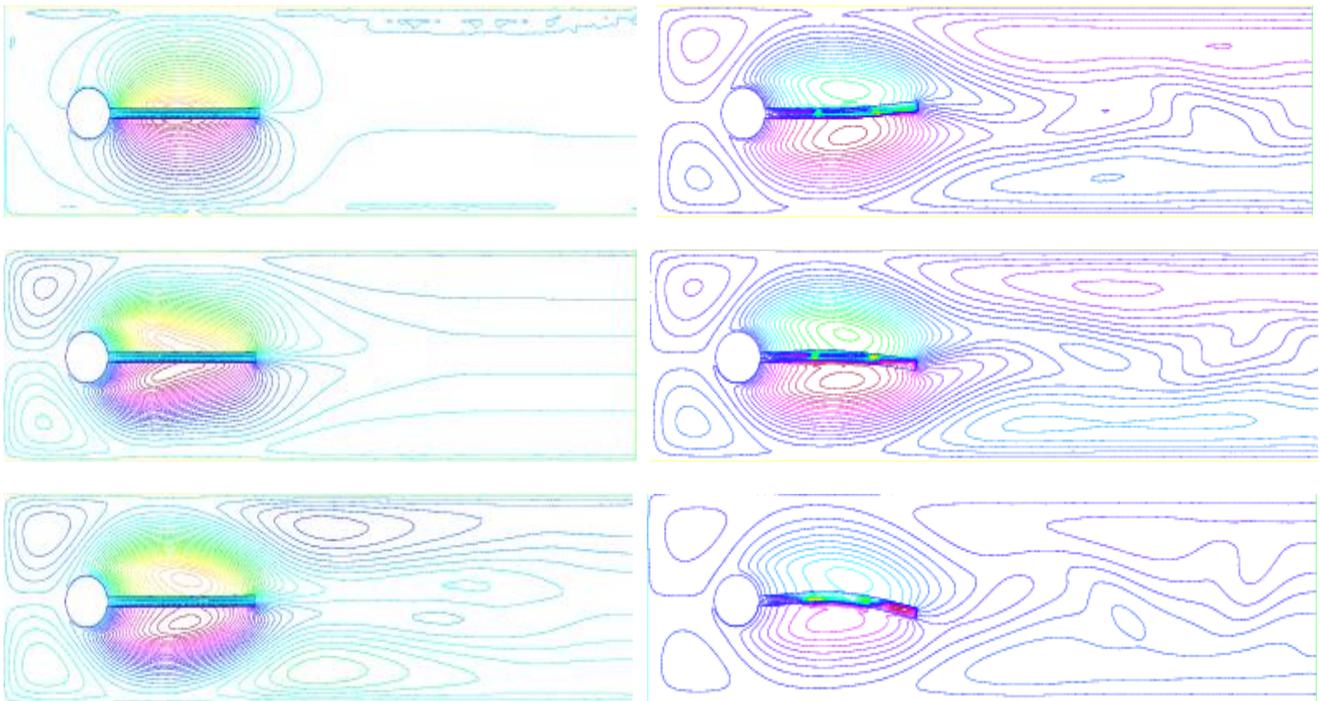

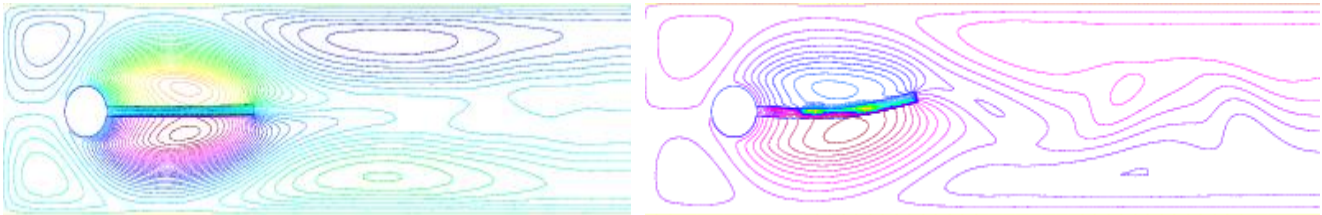

***Figure 7***: $w$- *contour plots* (*column* $1: t = 0.005, 0.2, 0.5, 0.75$ & *column* $2: t = 2.0, 2.1, 2.95, 3.05$ ).

## CONCLUSION

In this study the effects of micro-viscosity parmaters of Cosserat fluid at micro-structural level has been analyzed by presenting a non-classical Cosserat fluid-structure interaction CFSI problem in a monolithic Eulerian formulation. The finite element method and semi-implicit scheme have been used for discretizing space and time domains. The proposed formulation is implemented using FreeFem++. The obtained results suggest that the amplitude of oscillations and micro-rotaional viscosity $\mu_r$ varies inversely, which specifies that fluctuation of the control poin $A$ almost vanishes by dominating micro-rotational viscosity $\mu_r$ on classical viscosity $\mu$ of the fluid. Moreover, micro-viscosity parameter $\lambda_1$ exhibit a significant effect on micro-rotation filed as compare to velocity field. Desired accuracy in results can be greatly avhieved by using a larger mesh on enhanced computational resources.

**Author's Contribuitions:** Conceptualization (Muhammad Sabeel Khan); Writing (Muhammad Sabeel Khan, Nazim Hussain Hajano and Lisheng Liu); Methodology (Muhammad Sabeel Khan); Supervision (Muhammad Sabeel Khan).